\documentclass [12 pt, twoside]{amsart}
\usepackage{enumerate}

 \usepackage{amsmath,amsfonts,amsthm,amssymb}
 \usepackage[all]{xy}

\newtheorem{thm}{Theorem}[section]
\newtheorem{prop}{Proposition}[section] 
\newtheorem{lem}{Lemma}[section] 
\newtheorem{cor}{Corollary}[section]

\newcommand{\bbb}[1]{\mbox{\boldmath$#1$}}

\theoremstyle{definition}
\newtheorem*{Proof}{Proof}

\newcommand{\ld} {{\ldots}}
\newcommand{\sm} {{\smallsetminus}}
\newcommand{\thi} {{\theta}}
\newcommand{\de} {{\delta}}
\newcommand{\De} {{\varDelta}}
\newcommand{\si} {{\sigma}}
\newcommand{\Si} {{\varSigma}}
\newcommand{\la} {{\lambda}}
\newcommand{\La} {{\varLambda}}

\newcommand{\e} {{\varepsilon}}

\newcommand{\mi} {{\mu}}

\newcommand{\dis}{\displaystyle}
\newcommand{\ssum}{\sum\limits}
\newcommand{\dlim}{\displaystyle\lim}
\newcommand{\ct}{{\mathcal{T}}}
\newcommand{\ca}{{\mathcal{A}}}

\newcommand{\cp}{{\mathcal{P}}}
\newcommand{\ch}{{\mathcal{H}}}

\newcommand{\cb}{{\mathcal{B}}}

\newcommand{\ra}{{\rightarrow}}
\newcommand{\fa}{{\forall}}

\newcommand{\qb}{$\quad\blacksquare$}
\def\1{\it1\hspace*{-0.150cm}{\footnotesize{I}}}

\def\C{{\mathbb{C}}}
\def\Q{{\mathbb{Q}}}
\def\N{{\mathbb{N}}}

\numberwithin{equation}{section}

\begin{document}

\title[ Common hypercyclic functions ]{ Common hypercyclic functions for translation operators with large gaps II}

\author[Nikos Tsirivas]{Nikos Tsirivas}
\address{Department of Mathematics and Applied Mathematics, University of Crete, GR-700 13 Heraklion, Crete, Greece}
\email{tsirivas@uoc.gr}
\thanks{The research project is implemented within the framework of the Action "Supporting Postdoctoral Researchers" of
the Operational Program "Education and Lifelong Learning" (Action's Beneficiary: General Secretariat for
Research and Technology), and is co-financed by the European Social Fund (ESF) and the Greek State.}

\subjclass[2010]{47A16}

\date{}

\keywords{Hypercyclic operator, Common hypercyclic vectors, Translation operator}

\begin{abstract}
Let $\ch(\C)$ be the set of entire functions endowed with the topology of local uniform convergence. Fix a
sequence of non-zero complex numbers $(\la_n)$, $|\la_n|\to +\infty$, which satisfies the following property:
for every $M>0$ there exists a subsequence $(\mi_n)$ of $(\la_n)$ such that \smallskip

(i) $|\mi_{n+1}|-|\mi_n|>M$ for every $n=1,2,\ld$ and  \smallskip

(ii) $\limsup_{n\to +\infty}|\mi_n|\ssum^{+\infty}_{k=n}\dfrac{1}{|\mi_k|}=+\infty$ .\smallskip

We prove that there exists an entire function $f$ such that for every non-zero complex number $a$ with $|a|=1$
the set $\{ f(z+\la_na):n=1,2,\ldots \}$ is dense in $\ch(\C)$. Actually, the set of entire functions with the
previous property is residual in $\ch(\C)$. This largely extends Theorem 2.2 in \cite{7}.
\end{abstract}

\maketitle

\section{Introduction}
This paper is the third one in a series of papers, \cite{Tsi}, \cite{Tsi2}, which aims to a better understanding
of the phenomenon of common hypercyclic vectors for uncountable many hypercyclic operators of translation type.
The notion of hypercyclicity has been studied intensively the last twenty years and there is by now a well
developed theory on this subject, see for instance the two recent books \cite{2}, \cite{10}. Let us recall the
relevant definitions. A sequence $(T_n)$ of continuous and linear operators acting on a real or complex
topological vector space $X$ is called hypercyclic, if there exists a vector $x\in X$ so that the set $\{
T_n(x): n=1,2,\ldots \}$ is dense in $X$; in this case $x$ is called hypercyclic for $(T_n)$ and the symbol
$HC(\{ T_n \})$ stands for the set of all hypercyclic vectors for $(T_n)$. In the previous definition if the
sequence $(T_n)$ comes from the iterates of a single operator $T$, i.e. $T_n=T^n$, $n=1,2,\ldots $ then $T$ is
called hypercyclic, $x$ is called hypercyclic for $T$ and $HC(T)$ denotes the set of hypercyclic vectors for
$T$. Observe that, in the above situation the topological vector space $X$ is necessarily separable. For several
examples of hypercyclic operators including many classical operators, such as weighted shifts, differential
operators, adjoints of multipliers and so on, we refer to \cite{2}, \cite{10}.

Our interest here lies on a particular operator, the translation operator acting on the space $\ch(\C)$ of
entire functions endowed with the topology $\ct_u$ of local uniform convergence. Let $a$ be a non-zero number of
the complex plane $\mathbb{C}$. For obvious reasons, the operator $T_a:\ch(\C)\to \ch(\C)$ defined by
$T_a(f)(z):=f(z+a)$, for $f\in\ch(\C)$, is called translation. A classical result due to Birkhoff \cite{4} says
that $T_1$ is hypercyclic. Actually for every $a\in \mathbb{C}\setminus \{ 0\}$ the translation operator $T_a$
is hypercyclic. This means, that there exists an entire function $f$ whose positive integer translates
approximate every entire function, i.e. the set $\{ f(z+n): n=1,2,\ldots \}$ is dense in $(\ch(\C),\ct_u)$. As
an easy application of Baire's category theorem we have the following dichotomy: if $X$ is a separable
topological vector space and $T$ is a linear and continuous operator on $X$ then either $HC(T)=\emptyset$ or
$HC(T)$ is $G_{\de}$ and dense in $X$, see \cite{2}, \cite{10}. Recall that a subset $A$ of $X$ is called
$G_{\de}$ if it can be written as countable intersection of open sets. Therefore, for every $a\in
\mathbb{C}\setminus \{ 0\}$ the set $HC(T_a)$ is is $G_{\de}$ and dense in $(\ch(\C),\ct_u)$ and as an immediate
consequence of Baire's category theorem we have that the set $\bigcap_{n=1}^{+\infty} HC(T_{a_n})$ is non-empty
for every sequence $(a_n)$ of non-zero complex numbers. Our point of departure is the following extension of
Birkhoff's theorem due to Costakis and Sambarino \cite{6}: the set $\bigcap\limits_{a\in \mathbb{C}\setminus \{
0\} }HC(T_{a})$ is residual in $(\ch(\C),\ct_u)$, hence non-empty. The difficulty of proving such a result is,
of course, the uncountable range of $a$. Subsequently Costakis, in an attempt to generalize the previous result
established the following.

\begin{thm} \cite{7} \label{tCos}
Let $(\la_n)$ be a sequence of non-zero complex numbers with $|\la_n| \to +\infty$ which also satisfies the
following condition:

for every $M>0$ there exists a subsequence $(\mi_n)$ of $(\la_n)$ such that \smallskip

(i) $|\mi_{n+1}|-|\mi_n|>M$ for every $n=1,2,\ld$ and  \smallskip

(ii) $\ssum^{+\infty}_{n=1}\dfrac{1}{|\mi_n|}=+\infty$.\smallskip

Then the set $\bigcap\limits_{|a|=1}HC(\{T_{\la_na}\})$ is a $G_\de$ and dense in $(\ch(\C),\ct_u)$.
\end{thm}

The purpose of the present work is to extend the above theorem by allowing a wider class of sequences $(\la_n)$.
However, it is necessary to impose certain restrictions on $(\la_n)$ so that the conclusion of the above theorem
holds, as the following result from \cite{8} shows: if $(\la_n)$ is a sequence of non-zero complex numbers with
$\liminf_n\frac{|\la_{n+1}|}{|\la_n|}>2$ then $\bigcap\limits_{|a|=1}HC(\{T_{\la_na}\})=\emptyset$.

In all this work we fix a sequence $(\la_n)$ of non-zero complex numbers that tends to infinity and, in
addition, $(\la_n)$ satisfies the following condition which we call it condition $(C)$ from now on.\\

For every $M>0$ there exists a subsequence $(\mi_n)$ of $(\la_n)$ such that \smallskip

(i) $|\mi_{n+1}|-|\mi_n|>M$ for every $n=1,2,\ld$ and  \smallskip

(ii) $\dis\limsup_{n\ra+\infty}\Big(|\mi_n|\cdot\ssum^{+\infty}_{k=n}\dfrac{1}{|\mi_k|}\Big)=+\infty$.
\smallskip

We denote by $C_{r}:=\{z\in\C|\,|z|=r\}$ the circle with center $0$ and radius $r$. Our main task in this paper
is to prove the following
\begin{thm}\label{thm2.1}
Fix a sequence of non-zero complex numbers $\La=(\la_n)$ that tends to infinity and satisfies the above
condition $(C)$.\ Then for every $r\in (0,+\infty)$ the set $\bigcap\limits_{a\in C_{r}}HC(\{T_{\la_na}\})$ is a
$G_\de$ and dense in $(\ch(\C),\ct_u)$. In particular, $$\bigcap\limits_{a\in C_{r}}HC(\{T_{\la_na}\})\neq
\emptyset.$$
\end{thm}
It is clear that property $(ii)$ of condition $(C)$ in Theorem \ref{thm2.1} relaxes the corresponding condition
in Theorem \ref{tCos}. It is important to mention that in \cite{Tsi2} we obtain a full strength of the
conclusion of Theorem \ref{tCos}, namely we show that under the assumptions of Theorem \ref{tCos} the set
$\bigcap\limits_{a\in \mathbb{C}\setminus \{ 0\}}HC(\{T_{\la_na}\})$ is $G_\de$ and dense in $(\ch(\C),\ct_u)$.
On the other hand, relaxing condition $(ii)$ of Theorem \ref{tCos} as above, the price we pay, at least for now,
is the "thin", but still uncountable, range of $a$ in the conclusion of Theorem \ref{thm2.1}. A further
connection of the present work with our main result from \cite{Tsi} will be discussed in Section 5.

There are several recent results concerning either the existence or the non-existence of common hypercyclic
vectors for uncountable families of operators, see for instance, \cite{AbGo}, \cite{Ba1}-\cite{Ber},
\cite{ChaSa1}-\cite{GaPa}, \cite{10}, \cite{LeMu}, \cite{Math}, \cite{San}-\cite{Tsi2}.

The paper is organized as follows. Sections 2-4 contain the proof of Theorem \ref{thm2.1}. In the last section,
Section 5, we connect our work with the main results from \cite{Tsi}, \cite{Tsi2}.

\section{Three Basic Lemmas}
Let us now describe the main steps for the proof of Theorem \ref{thm2.1}. Defining the arcs
\[
A_k:=\bigg\{a\in\C|\exists\;t\in\bigg[\frac{k}{4},\frac{k+1}{4}\bigg] \ \ \text{such that} \ \ a=r_0e^{2\pi
it}\bigg\}, \ \ k=0,1,2,3
\]
and using Baire's category theorem we easily see that Theorem \ref{thm2.1} reduces to the following
\begin{prop}\label{prop2.1}
Fix a sequence $(\lambda_n)$ of non-zero complex numbers that tends to infinity and satisfies the above
condition $(C)$. Fix three real numbers $r_0,\thi_0,\thi_T$ such that $r_0\in(0,+\infty)$,
$0\le\thi_0<\thi_T\le1$, $\thi_T-\thi_0=\dfrac{1}{4}$ and consider the arc $A$ defined by
\[
A=\{a\in\C|\,\text{there exists}\;t\in[\thi_0,\thi_T] \; \text{such that}\; a=r_0e^{2\pi it}\}.
\]
Then $\bigcap\limits_{a\in A}HC(\{T_{\la_na}\})$ is a $G_\de$ and dense subset of $(\ch(\C),\ct_u)$.
\end{prop}
For the proof of Proposition \ref{prop2.1} we introduce some notation which will be carried out throughout this
paper. Let $(p_j)$, $j=1,2,\ld$ be a dense sequence of $(\ch(\C),\ct_u)$, (for instance, all the polynomials in
one complex variable with coefficients in $\Q+i\Q$). For every $m,j,s,k\in\N$ we consider the set $$E(m,j,s,k):=
\!\Big\{f\!\in\!\ch(\C)\,\, | \forall a\!\in\! A \,\,\,\exists  n\!\in\!\N, n\!\le\! m : \dis\sup_{|z|\le
k}|f(z\!+\!\la_na)\!-\!p_j(z)|\!<\!\frac{1}{s}\Big\}.$$
Clearly, Baire's category theorem and the three lemmas stated below imply Proposition \ref{prop2.1}.
\begin{lem}\label{lem2.1}
\[
\bigcap_{a\in A}HC(\{ T_{\la_na}\} )=\bigcap^\infty_{j=1}\bigcap^\infty_{s=1}\bigcap^\infty_{k=1}
\bigcup^\infty_{m=1}E(m,j,s,k).
\]
\end{lem}
\begin{lem}\label{lem2.2}
For every $m,\!j,\!s,\!k\!\in\!\N$ the set $E(m,\!j,\!s,\!k)$ is open in $(\ch(\C),\ct_u)$.
\end{lem}
\begin{lem}\label{lem2.3}
For every $j,\!s,\!k\in\N$ the set $\bigcup\limits^\infty_{m=1}E(m,j,s,k)$ is dense in $(\ch(\C),\ct_u)$.
\end{lem}

The proof of Lemma \ref{lem2.1} is in \cite{Tsi}. The proof of Lemma \ref{lem2.2} is similar to that in Lemma 9
of \cite{6} and it is omitted. It remains to prove Lemma \ref{lem2.3}. This will be done in Sections 3 and 4.
\section{Construction of the partition and the disks}\label{sec3}
\noindent

Let $\La=(\la_n)$ be a sequence of non-zero complex numbers such that $\la_n\ra\infty$ as $n\ra+\infty$ and we
further assume that $\La$ satisfies condition $(C)$. For the sequel we fix four positive numbers
$c_1,c_2,c_3,c_4$ such that $c_1>1$, $c_2\in(0,1)$, $c_3>1$, $c_4>1$, where $c_3:=\dfrac{c_4}{r_0c_2}$,
$c_1:=4(c_3+1)$. Let three real numbers $\thi_0,\thi_T,r_0$ be as in Proposition \ref{prop2.1}. After the
definition of the above numbers we fix a subsequence $(\mi_n)$ of $(\la_n)$ such that:
\[
|\mi_n|,|\mi_{n+1}|-|\mi_n|>c_1 \ \ \text{for every} \ \ n=1,2,\ld \ \ \text{and} \ \
\lim_{n\ra+\infty}\bigg(|\mi_n|\cdot\sum^{+\infty}_{k=n} \frac{1}{|\mi_k|}\bigg)=+\infty.
\]
The condition $\dis\lim_{n\ra+\infty}\Big(|\mi_n|\cdot\ssum^{+\infty}_{k=n}\dfrac{1}{|\mi_k|}\Big)=+\infty$
implies that there exists a positive integer $m_0$ such that for every $m\geq m_0$
\[
\sum_{k=m}^{+\infty}\frac{1}{|\mi_k|}>c_3\cdot\frac{1}{|\mi_m|}.
\]
Throughout Section 3 the positive integer $m_0$ will appear frequently and it is fixed from now on.

\subsection{Step 1. Partitions of the interval $\bbb{[\thi_0,\thi_T]}$} \label{subsec3.1}
\noindent

For every sufficiently large positive integer $m$ (actually $m\geq m_0$) we shall construct a corresponding
partition $\De_m$ of $[\thi_0,\thi_T]$. For every $m\geq m_0$ let $m_1(m)$ be the minimum positive integer such
that:
\begin{eqnarray}
\sum^{m_1(m)}_{k=m}\frac{1}{|\mi_k|}>c_3\cdot\frac{1}{|\mi_m|}.  \label{sec3eq1}
\end{eqnarray}
Obviously $m_1(m)\ge m+1$ for every $m=m_0, m_0+1\ld$, since $c_3>1$. Let $m$ be any positive integer with
$m\geq m_0$. We define the numbers $\thi_0^{(m)}:=\thi_0$, $\thi^{(m)}_1:=\thi^{(m)}_0+\dfrac{c_2}{|\mi_m|}$,
$\thi^{(m)}_2:=\thi^{(m)}_1+\dfrac{c_2}{|\mi_{m+1}|},\ld$,
$\thi^{(m)}_{m_1(m)-m+1}:=\thi^{(m)}_{m_1(m)-m}+\dfrac{c_2}{|\mi_{m_1(m)}|}$, or in a more compact form
\begin{eqnarray}
\thi^{(m)}_{n+1}:=\thi^{(m)}_n+\frac{c_2}{|\mi_{m+n}|} \ \ \text{for} \ \ n=0,1,\ld, \ \ m_1(m)-m.
\label{sec3eq2}
\end{eqnarray}
We denote $$\si_m:=\thi^{(m)}_{m_1(m)-m+1}-\thi_0.$$ Now let some positive integer $\nu>m_1(m)-m+1$, $\nu\in\N$.
Then there exists a unique pair $(k,j)\in\N^2$, where $j\in\{0,1,\ld,m_1(m)-m\}$ such that:
\[
\nu=k(m_1(m)-m+1)+j.
\]
Define
\[
\thi^{(m)}_\nu:=\thi^{(m)}_j+k\si_m, \ \ \nu> m_1(m)-m+1.
\]
It is obvious that $\dlim_{\nu\ra+\infty}\thi^{(m)}_\nu=+\infty$ and the sequence $(\thi^{(m)}_\nu)_\nu$ is
strictly increasing in respect to $\nu$. So there exists a maximum natural number $\nu_m\in\N$ such that
$\thi^{(m)}_{\nu_m}\le\thi_T$. We set $\De_m:=\{\thi^{(m)}_0,\thi^{(m)}_1,\ld,\thi^{(m)}_{\nu_m}\}$. It holds
$\nu_m\ge m_1(m)-m+1$ (see Lemma \ref{sec3.4lem1}).
\subsection{Step 2. Partitions of the arc $\bbb{\phi_{r_0}([\theta_0 ,\theta_T])=A}$} \label{subsec3.2}
\noindent

Consider the function $\phi :[\theta_0 ,\theta_T ]\times (0 ,+\infty )\to \mathbb{C}$ given by
$$\phi (t,r_0):=r_0e^{2\pi it},\,\,\, (t,r_0)\in [\theta_0 ,\theta_T ]\times (0 ,+\infty ) $$ and for $r_0> 0$ we
define the corresponding curve $\phi_{r_0}:[\theta_0 ,\theta_T ] \to \mathbb{C}$ by
$$\phi_{r_0}(t):=\phi (t,r_0), \,\,\, t\in [\theta_0 ,\theta_T ].$$
For any given positive integer $m\geq m_0$, $\phi_{r_0} (\Delta_m)$ is a partition of the arc
$\phi_{r_0}([\theta_0 ,\theta_T])$, where $m_0$ is the positive integer defined above in Step 1 and $\Delta_m$
is the partition of the interval $[\theta_0 ,\theta_T]$ constructed in Step 1. For every $m\in \mathbb{N}$ with
$m\geq m_0$ define
$$\cp_{m}:=\phi_{r_0}(\Delta_m)$$
which we call partition of the arc $\phi_{r_0}([\theta_0 ,\theta_T])=A$. %

\subsection{Step3. A lower bound for the stopping time} \label{subsec3.4}

\begin{lem} \label{sec3.4lem1}
Let $m\in\N$ with $m\geq m_0$. Then $\si_m=\thi^{(m)}_{m_1(m)-m+1}-\thi_0<1/4$. In particular $\nu_m\ge
m_1(m)-m+1$.
\end{lem}
\begin{Proof}
By the definition of the numbers $\thi^{(m)}_j$, $j=0,1,\ld,m_1(m)-m+1$ we have
\begin{eqnarray}
\thi^{(m)}_{m_1(m)-m+1}-\thi_0=c_2\cdot\sum^{m_1(m)}_{k=m}\frac{1}{|\mi_k|}.  \label{sec3eq2}
\end{eqnarray}
In order to bound the right hand term in the equality above, observe that
\begin{eqnarray} \label{sec3eq3}
\sum^{m_1(m)}_{k=m}\frac{1}{|\mi_k|}\le
c_3\cdot\frac{1}{|\mi_m|}+\frac{1}{|\mi_{m_1(m)}|}<(c_3+1)\frac{1}{|\mi_m|},
\end{eqnarray}
which follows from the definition of the number $m_1(m)$.  Since $c_1=4(c_3+1)$ and
$|\mi_m|>c_1=4(c_3+1)>4c_2(c_3+1)$ , recall that $c_2\in(0,1)$,
\begin{eqnarray}
\frac{c_3+1}{|\mi_m|}<\frac{1}{4c_2}.  \label{sec3eq4}
\end{eqnarray}
Thus by (\ref{sec3eq2}), (\ref{sec3eq3}) and (\ref{sec3eq4}), $\si_m<\dfrac{1}{4}$ and the proof is complete.
\qb
\end{Proof}

\noindent

\subsection{Step 4. Construction of the disks} \label{subsec4.1}
\noindent

Our task in this subsection is to assign to each point $w$ of the partition $\cp_m$ for $m\geq m_0$ a suitable
closed disk with center $w\mi(w)$ and radius $c_4$ (the radius will be the same for every member of the family
of the disks), where $\mi(w)$ will be chosen from the sequence $(\mi_n)$. We shall see that, the construction of
the partition $\cp_m$ ensures on the one hand that the points of the partition are close enough to each other on
the arc $A$ and on the other hand that the disks centered at these points with fixed radius $c_4$ are pairwise
disjoint.

We set
\[
\cb:=\{z\in\C/|z|\le c_4\}
\]
and fix any positive integer $m$ with $m\geq m_0$. Let $w$ be an arbitrary point in $\cp_m$. There exists unique
$n\in\{0,1,\ld,\nu_m\}$ such that $w=r_0e^{2\pi i\thi^{(m)}_n}$. Now there exists unique $k\in\N$, $k\ge1$ and
$j\in\{0,1,\ld,m_1(m)-m\}$ such that $n=k(m_1(m)-m+1)+j$ and define $$\mi(w):=\mi_{m +j}.$$ Thus we assign, in a
unique way, a term of the sequence $(\mi_n)$ to every point of $\cp_m$ and to conclude our construction we
introduce the notation
$$\cb_w:=\cb+w\mi(w).$$
The desired disks are the disks $\cb$ and $\cb_w$, $w\in\cp_m$. Denote by
$$\mathfrak{D}_m:=\{\cb\}\cup\{\cb_w:w\in\cp_m\} $$ the collection of the above disks.
\subsection{Step 5. The disks are pairwise disjoint}\label{subsec4.2}
\begin{lem}\label{lem4.1}
Let $m\in \mathbb{N}$ with $m\geq m_0$. Then $\cb\cap\cb_w=\emptyset$ for every $w\in\cp_m$.
\end{lem}
\begin{Proof}
We have $c_3=\dfrac{c_4}{r_0c_2}>\dfrac{c_4}{r_0}$, since $c_2\in(0,1)$, and taking into account
$c_1=4(c_3+1)>2c_3$ we get
\begin{eqnarray}
c_1>\frac{2c_4}{r_0}.  \label{sec3eq5}
\end{eqnarray}
Let $w\in\cp_m$. The closed disks $\cb$, $\cb_w$ are centered at $0$, $w\mi(w)$ respectively and they have the
same radius $c_4$. Hence, we have to show that $|w\mi(w)|>2c_4$. Since $|w|\ge r_0$, it suffices to prove that
$|\mi(w)|>\dfrac{2c_4}{r_0}$. Observe now that, by the definition of $\mi(w)$ in the previous subsection,
\begin{eqnarray}
\mi(w)=\mi_n  \label{sec3eq6}
\end{eqnarray}
for some positive integer $n\in\N$. As a consequence of the definition of the sequence $(\mi_n)$ we have
\begin{eqnarray}
|\mi_n|>c_1 \ \ \text{for every} \ \ n\in\N.  \label{sec3eq7}
\end{eqnarray}
By (\ref{sec3eq5}), (\ref{sec3eq6}) and (\ref{sec3eq7}) we conclude that $|\mi(w)|>2\dfrac{c_4}{r_0}$ and this
finishes the proof of the lemma. \qb
\end{Proof}
\begin{lem}\label{lem4.2}
Let $m\in \mathbb{N}$ with $m\geq m_0$ and $w_1,w_2\in\cp_m$ with $w_1\neq w_2$. Then
$\cb_{w_1}\cap\cb_{w_2}=\emptyset$.
\end{lem}
\begin{Proof}
We distinguish two cases:\smallskip

(i) $|\mi(w_1)|<|\mi(w_2)|$. \smallskip

Our hypothesis implies
\begin{align*}
|w_2\mi(w_2)-w_1\mi(w_1)|&\ge\big||w_2\mi(w_2)|-|w_1\mi(w_1)|\big|\\
&=|w_1| \cdot(|\mi(w_2)|-|\mi(w_1)|)\ge r_0\cdot c_1>2c_4,
\end{align*}
thus $\cb_{w_1}\cap\cb_{w_2}=\emptyset$.\smallskip

(ii) $|\mi(w_1)|=|\mi(w_2)|$. \smallskip

By the definition of the partition $\cp_m$ we have $w_1=r_0\cdot e^{2\pi i\thi^{(m)}_{n_1}}$, $w_2=r_0\cdot
e^{2\pi i\thi^{(m)}_{n_2}}$ for some $n_1,n_2\in\{0,1,\ld,\nu_m\}$ and $n_1\neq n_2$ because $w_1\neq w_2$.
Without loss of generality we suppose that $n_1<n_2$. Now there exists a unique pair $(k_1,j_1)$, where
$k_1\in\N$, $j_1\in\{0,1,\ld,m_1(m)-m\}$ and a unique pair $(k_2,j_2)$ where $k_2\in\N$ and
$j_2\in\{0,1,\ld,m_1(m)-m\}$ such that
\begin{eqnarray}
n_1=k_1(m_1(m)-m+1)+j_1 \ \  \label{sec3eq8}
\end{eqnarray}
and
\begin{eqnarray}
n_2=k_2(m_1(m)-m+1)+j_2.  \label{sec3eq9}
\end{eqnarray}
By definition we have $\mi(w_1)=\mi_{m+j_1}$, $\mi(w_2)=\mi_{m+j_2}$ and the hypothesis yields
\[
|\mi(w_1)|=|\mi(w_2)|\Leftrightarrow\mi(w_1)=\mi(w_2).
\]
So we have $j_1=j_2=j_0$.
Thus
\[
\thi^{(m)}_{n_1}=\thi^{(m)}_{j_0}+k_1\si_m ,
\]
\[
\thi^{(m)}_{n_2}=\thi^{(m)}_{j_0}+k_2\si_m
\]
and
\begin{eqnarray}
\thi^{(m)}_{n_2}-\thi^{(m)}_{n_1}=(k_2-k_1)\si_m.  \label{sec3eq10}
\end{eqnarray}
By (\ref{sec3eq8}), (\ref{sec3eq9}) and the fact that $n_1<n_2$ and $j_1=j_2$ we have $k_1<k_2\Rightarrow k_2\ge
k_1+1$. Using now (\ref{sec3eq10}) it follows that
\begin{eqnarray}
\thi^{(m)}_{n_2}-\thi^{(m)}_{n_1}\ge\si_m>0.  \label{sec3eq11}
\end{eqnarray}
A lower bound for the quantity $|w_2\mi(w_2)-w_1\mi(w_1)|$ is:
\begin{align}
|w_2\mi(w_2)-w_1\mi(w_1)|&=|\mi(w_1)|\cdot|w_1-w_2|\ge|\mi_m|\cdot|w_1-w_2| \nonumber\\
&=|\mi_m|\cdot|r_0\cdot e^{2\pi i\thi^{(m)}_{n_2}}-r_0e^{2\pi i\thi^{(m)}_{n_1}}| \nonumber \\
&=r_0|\mi_m|\cdot|e^{2\pi i\thi^{(m)}_{n_2}}-e^{2\pi i\thi^{(m)}_{n_1}}| \nonumber \\
&=r_0|\mi_m|\cdot 2\sin(\pi(\thi^{(m)}_{n_2}-\thi^{(m)}_{n_1})). \nonumber\\
\label{sec3eq12}
\end{align}
Consider Jordan's inequality $$\sin x>\dfrac{2}{\pi}x, \,\,\, x\in\Big(0,\dfrac{\pi}{2}\Big).$$ We have
\[
0<\thi^{(m)}_{n_2}-\thi^{(m)}_{n_1}<\frac{1}{4}\Rightarrow0<\pi(\thi^{(m)}_{n_2}-
\thi^{(m)}_{n_1})<\frac{\pi}{4}.
\]
So, applying Jordan's inequality for $x=\pi(\thi^{(m)}_{n_2}-\thi^{(m)}_{n_1})$ we derive
\begin{eqnarray} \label{sec3eq13}
\sin(\pi(\thi^{(m)}_{n_2}-\thi^{(m)}_{n_1}))
>2(\thi^{(m)}_{n_2}-\thi^{(m)}_{n_1}).
\end{eqnarray}
Now, inequalities (\ref{sec3eq11}), (\ref{sec3eq12}) and (\ref{sec3eq13}) imply
\begin{eqnarray}
|w_2\mi(w_2)-w_1\mi(w_1)|>4r_0|\mi_m|\cdot\si_m.  \label{sec3eq14}
\end{eqnarray}
By the definition of the number $\si_m$ and relation (\ref{sec3eq2}) of Lemma \ref{sec3.4lem1} we obtain
\[
\si_m=c_2\cdot\sum^{m_1(m)}_{k=m}\frac{1}{|\mi_k|}.
\]
The last equality, inequality (\ref{sec3eq14}) and the definition of the number $m_1(m)$ give
\begin{align*}
|w_2\mi(w_2)-w_1\mi(w_1)|&>4r_0|\mi_m|\cdot c_2\sum^{m_1(m)}_{k=m}
\frac{1}{|\mi_k|} \nonumber \\
&>4r_0|\mi_m|\cdot c_2\frac{c_3}{|\mi_m|}=4r_0c_2c_3
\end{align*}
and the properties of our fixed numbers imply $4r_0c_2c_3>2c_4$. It follows that
$\cb_{w_1}\cap\cb_{w_2}=\emptyset$ and the proof of this lemma is complete. \qb
\end{Proof}

By Lemmas \ref{lem4.1}, \ref{lem4.2} we conclude the following
\begin{cor}\label{cor4.3}
For every positive integer $m$ with $m\geq m_0$ the family $\mathfrak{D}_m:=\{\cb\}\cup\{\cb_w:w\in\cp_m\}$
consists of pairwise disjoint disks. \qb
\end{cor}
\section{Proof of Lemma \ref{lem2.3}}\label{sec4}
\noindent

Let $j_1,s_1,k_1\in\N$ be fixed. We will prove that the set $\bigcup\limits^\infty_{m=1}E(m,j_1,s_1,k_1)$ is
dense in $(\ch(\C),\ct_u)$. For simplicity we write $p_{j_1}=p$. Consider fixed $g\in\ch(\C)$, a compact set
$C\subseteq\C$ and $\e_0>0$. We seek $f\in\ch(\C)$ and a positive integer $m_1$ such that
\begin{equation}\label{sec4eq1}
f\in E(m_1,j_1,s_1,k_1)
\end{equation}
and
\begin{equation} \label{sec4eq2}
\sup_{z\in C}|f(z)-g(z)|<\e_0.
\end{equation}
Fix $R_1>0$ sufficiently large so that
$$
C\cup\{z\in\C|\,|z|\le k_1\}\subset\{z\in\C|\,|z|\le R_1\}.
$$
and then choose $0<\de_0<1$ such that
\begin{equation} \label{sec4eq3}
\textrm{if} \,\,\, |z|\le R_1 \,\,\, \textrm{and} \,\,\, |z-w|<\de_0, \ \ w\in\C, \,\,\, \textrm{then} \,\,\,
|p(z)-p(w)|<\frac{1}{2s_1}.
\end{equation}
We set $$\cb:=\{z\in\C|\,|z|\le R_1+\de_0\},$$
\[
c_4:=R_1+\de_0, \quad c_2:=\frac{\de_0}{2(2r_0\pi+1)},
\]
\[
c_3=\frac{c_4}{r_0c_2}=\frac{R_1+\de_0}{r_0\dfrac{\de_0}{2(2r_0\pi+1)}}=
\frac{2(R_1+\de_0)(2r_0\pi+1)}{r_0\de_0}
\]
\[
c_1=4(c_3+1)=4\cdot\bigg(\frac{2(R_1+\de_0)(2r_0\pi+1)}{r_0\de_0}+1\bigg).
\]
After the definition of the above numbers we choose a subsequence $(\mi_n)$ of $(\la_n)$ such that\smallskip

(i) $|\mi_n|,|\mi_{n+1}|-|\mi_n|>c_1$ for $n=1,2,\ld$  \bigskip

and

(ii) $\dis\lim_{n\ra+\infty}\bigg(|\mi_n|\cdot\ssum^{+\infty}_{k=n}\dfrac{1}{|\mi_k|}\bigg) =+\infty$. \bigskip

By condition (ii), there exists a positive integer $m_0$ such that for every $m\geq m_0$
\[
\sum_{k=m}^{+\infty}\frac{1}{|\mi_k|}>c_3\cdot\frac{1}{|\mi_m|}.
\]
After that, on the basis of the fixed numbers $r_0,\thi_0,\thi_T,c_1,c_2,c_3,c_4$ and $m_0$, for every positive
integer $m$ with $m\geq m_0$ we define the corresponding partition $\cp_m$ (see section 3). From now on till the
end of the proof of the lemma we fix a positive integer $m\geq m_0$ with its corresponding partition $\cp_m$.
For simplicity we write
$$\cp:=\cp_m.$$
Then, we define the set $L$ as follows:
$$L:=\cb\cup \left( \bigcup_{w\in \cp}\cb_w \right) ,$$
where the discs $\cb_w$, $w\in \cp$ are constructed in Section 3. By Corollary \ref{cor4.3}, the family
$\mathfrak{D}:=\{ \cb \} \cup \{ \cb_w : w\in \cp \}$ consists of pairwise disjoint disks. Therefore the compact
set $L$ has connected complement. This property is needed in order to apply Mergelyan's theorem. We now define
the function $h$ on the compact set $L$, $h:L\ra\C$ by
\[
h(z)=\left\{\begin{array}{cc}
              g(z), & z\in \cb \\
              p(z-w\lambda (w)), & z\in \cb_w ,w\in \cp.
            \end{array}\right.
\]
By Mergelyan's theorem \cite{12} there exists an entire function $f$ (in fact a polynomial) such that
\begin{eqnarray}
\sup_{z\in L}|f(z)-h(z)|<\min\bigg\{\frac{1}{2s_1},\e_0\bigg\}.  \label{sec4eq4}
\end{eqnarray}
By the definition of $h$, (\ref{sec4eq4}) and the definitions of sets $C$ and $\cb$, where $C\subseteq\cb$, it
follows that
\begin{eqnarray*}
\sup_{z\in C}|f(z)-g(z)|\le\sup_{z\in\cb}|f(z)-g(z)|\le\sup_{z\in L}|f(z)-h(z)|<\e_0,
\end{eqnarray*}
which implies the desired inequality (\ref{sec4eq2}).

It remains to show (\ref{sec4eq1}). Let some $a\in A$. There exists $\thi\in[\thi_0,\thi_T]$ such that
$a=r_0e^{2\pi i\thi}$. Now there exists unique $\rho\in\{0,1,\ld,\nu_m-1\}$ such that:
\[
\text{either}\,\,\, \thi^{(m)}_\rho\le\thi<\thi^{(m)}_{\rho+1} \ \ \text{or} \ \
\thi^{(m)}_{\nu_m}\le\thi\le\thi_T.
\]
and then define
\[
\thi_1:=\thi^{(m)}_\rho \,\,\,\textrm{and}\,\,\, \thi_2:=\thi^{(m)}_{\rho+1} \,\,\, \text{if}\,\,\,
\thi^{(m)}_\rho\le\thi<\thi^{(m)}_{\rho+1},
\]
\[
\thi_1:=\thi^{(m)}_{\nu_m} \,\,\,\textrm{and}\,\,\,  \thi_2:=\thi_T \,\,\, \text{if}\,\,\,
\thi^{(m)}_{\nu_m}\le\thi\le\thi_T.
\]
For the above, recall the definitions of $\nu_m$ and $\thi^{(m)}_\rho$, $\rho\in\{0,1,\ld,\nu_m-1\}$ from
Section 3. Set $$w_0:=r_0\cdot e^{2\pi i\thi_1}\in\cp.$$ We shall prove that for every $z\in\C$, $|z|\le R_1$ we
have $z+a\mi(w_0)\in\cb_{w_0}$. Recall that $\cb_{w_0}:=\cb+w_0\mi(w_0)=\overline{D}(w_0\mi(w_0),R_1+\de_0)$. It
suffices to prove that
\begin{eqnarray}
|(z+a\mi(w_0))-w_0\mi(w_0)|<R_1+\de_0 \ \ \text{for} \ \ |z|\le R_1.  \label{sec4eq5}
\end{eqnarray}
For $|z|\le R_1$ we have:
\begin{align}
|(z+a\mi(w_0))-w_0\mi(w_0)|&\le R_1+|\mi(w_0)|\,|a-w_0|\nonumber \\
&=R_1+|\mi(w_0)|\cdot| r_0\cdot e^{2\pi i\thi}-r_0e^{2\pi i\thi_1}|.  \label{sec4eq6}
\end{align}
By (\ref{sec4eq5}) and (\ref{sec4eq6}) it suffices to prove
\begin{eqnarray}
|\mi(w_0)|\cdot|r_0e^{2\pi i\thi}-r_0e^{2\pi i\thi_1}|<\de_0.  \label{sec4eq7}
\end{eqnarray}
We have:
\begin{align*}
|r_0e^{2\pi i\thi}-r_0e^{2\pi i\thi_1}|&=r_0|e^{2\pi i\thi}-e^{2\pi i\thi_1}| \nonumber \\
&=r_02\sin(\pi(\thi-\thi_1)) \nonumber\\
&\le r_02\sin(\pi(\thi_2-\thi_1)) \nonumber\\
&\le2r_0\pi(\thi_2-\thi_1)\nonumber\\
&\le2r_0\pi\dfrac{c_2}{|\mi(w_0)|} \nonumber\\
&=2r_0\pi\dfrac{\de_0}{2(2r_0\pi+1)}\cdot\dfrac{1}{|\mi(w_0)|}<\dfrac{\de} {2|\mi(w_0)|},
\end{align*}
which implies (\ref{sec4eq7}). So, for every $z\in\C$, $|z|\le R_1$
\begin{eqnarray}
z+a\mi(w_0)\in\cb_{w_0}.  \label{sec4eq8}
\end{eqnarray}
The definition of $h$ and (\ref{sec4eq8}) give that for every $z\in\C$, $|z|\le R_1$
\begin{eqnarray}
|f(z+a\mi(w_0))-p(z+\mi(w_0)(r_0e^{2\pi i\thi}-r_0e^{2\pi i\thi_1})|<\frac{1}{2s_1}.  \label{sec4eq9}
\end{eqnarray}
By (\ref{sec4eq3}) and (\ref{sec4eq7}) we get: for every $z\in\C$, $|z|\le R_1$
\begin{eqnarray}
|p(z+\mi(w_0)(r_0e^{2\pi i\thi}-r_0e^{2\pi i\thi_1}))-p(z)|<\frac{1}{2s_1}.  \label{sec4eq10}
\end{eqnarray}
The triangle inequality, for $z\in\C$, $|z|\le R_1$, gives
\begin{align}
|f(z+a\mi(w_0))-p(z)|\le&|f(z+a\mi(w_0))-p(z+\mi(w_0)(r_0e^{2\pi i\thi}-r_0
e^{2\pi i\thi_1}))| \nonumber\\
&+|p(z+\mi(w_0)(r_0e^{2\pi i\thi}-r_0e^{2\pi i\thi_1}))-p(z)|.  \label{sec4eq11}
\end{align}
By (\ref{sec4eq9}), (\ref{sec4eq10}), (\ref{sec4eq11}) and the fact that $k_1\le R_1$ we arrive at
\begin{eqnarray}
\sup_{|z|\le k_1}|f(z+a\mi(w_0))-p(z)|<\frac{1}{s_1}.  \label{sec4eq12}
\end{eqnarray}
Setting
\[
m_1:=\max\{n\in\N|\la_n=\mi(w), \ \ \text{for some} \ \ w\in\cp\},
\]
observing that the definition of $m_1$ is independent from $a\in A$ and in view of (\ref{sec4eq12}) we conclude
that for every $a\in A$ there exists some $n\in\N$ with $n\le m_1$ such that
\[
\sup_{|z|\le k_1}|f(z+a\la_n)-p(z)|<\frac{1}{s_1},
\]
where $f\in\ch(\C)$, since $f$ is a polynomial. This implies (\ref{sec4eq1}) and the proof of the lemma is
complete. \qb
\section{Examples  of sequences $\bbb{\La:=(\la_n)}$ satisfying the condition $(C)$}
\noindent

In this section we show that our main theorem is not covered by our recent results in \cite{Tsi}, \cite{Tsi2}.
We say that a sequence of non-zero complex numbers $(\la_n)$ with $\la_n \to \infty$ satisfies condition $(\Si)$
if:

for every $M>0$ there exists a subsequence $(\mi_n)$ of $(\la_n)$ such that \smallskip

(i) $|\mi_{n+1}|-|\mi_n|>M$ for every $n=1,2,\ld$ and  \smallskip

(ii) $\ssum^{+\infty}_{n=1}\dfrac{1}{|\mi_n|}=+\infty$.\smallskip

We introduce the following definitions.
\[
\mathcal{L}=\{\La=(\la_n)\in\C^\N|\la_n\neq0\;\fa\;n\in\N \ \ \text{and} \ \ \la_n\ra\infty\},
\]
\[
\ca_1:=\{\La=(\la_n)\in \mathcal{L}| \La \,\,\, \textrm{satisfies condition} \,\,\, (\Si) \},
\]
\[
\ca_2:=\{\La=(\la_n)\in \mathcal{L}| \La \,\,\, \textrm{satisfies condition} \,\,\, (C) \}.
\]

For $\La\in \mathcal{L}$, define
$$\cb(\La):=$$
$$
\Big\{a\in[0,+\infty]| a=\underset{n\ra+\infty}{\lim\sup}\Big|\dfrac{\mi_{n+1}}{\mi_n}\Big|
\,\,\, \textrm{for some subsequence}\,\,\,  (\mi_n) \,\,\, \textrm{of} \,\,\, \La \Big\}
$$
and
$$i(\La):=\inf\cb(\La).$$
Clearly,
$$i(\La) \in [1,+\infty] \,\,\, \textrm{for every}\,\,\,\La\in \mathcal{L}.$$

Our main results in \cite{Tsi}, \cite{Tsi2} are the following
\begin{thm} [\cite{Tsi}]\label{theo5.1}
If $\La:=(\la_n)\in \mathcal{L}$ and $i(\La)=1$ then $\bigcap\limits_{a\in \mathbb{C} \setminus \{ 0\}
}HC(\{T_{\la_na}\})$ is a $G_\de$ and dense subset of $(\ch(\C),\ct_u)$.
\end{thm}

\begin{thm} [\cite{Tsi2}]\label{theo5.2}
If $\La:=(\la_n)\in \ca_1$ then $\bigcap\limits_{a\in \mathbb{C} \setminus \{ 0\} }HC(\{T_{\la_na}\})$ is a
$G_\de$ and dense subset of $(\ch(\C),\ct_u)$.
\end{thm}

In view of Theorem \ref{theo5.1} and in order to completely characterize the sequences $\La:=(\la_n)\in
\mathcal{L}$ such that the set $\bigcap\limits_{a\in \mathbb{C} \setminus \{ 0\} }HC(\{T_{\la_na}\})$ is a
$G_\de$ and dense in $(\ch(\C),\ct_u)$ one has to deal with sequences $\La\in \mathcal{L}$ for which $i(\La)>1$.
This is one of the reasons we introduced the classes $\ca_1$, $\ca_2$. Indeed, it is established in \cite{Tsi2}
that the class $\ca_1$ contains sequences $\La\in \mathcal{L}$ with $i(\La)>1$. On the other hand, there exist
sequences $\La\in \mathcal{L}$ with $i(\La)=1$ and $\La\notin \ca_1$, see \cite{Tsi2}. Since,
$$\ca_1 \subset \ca_2$$ we conclude that the class $\ca_2$ contains sequences $\La\in \mathcal{L}$ with
$i(\La)>1$. The above inclusion is strict; for instance, the sequence $\la_n=n^2$, $n=1,2,\ld$, belongs to
$\ca_2$ but not in $\ca_1$. However $i((n^2))=1$, therefore for this sequence the conclusion of Theorem
\ref{theo5.1} holds and of course in this case the conclusion of Theorem \ref{thm2.1} is covered by the much
stronger Theorem \ref{theo5.1}. So the interest here is to show that there exist $\La\in\ca_2\sm A_1$ with
$i(\La)=M$ for some positive real number $M>1$, and this in turn shows that our main result, Theorem
\ref{thm2.1}, is not covered by Theorems \ref{theo5.1}, \ref{theo5.2}. This is the content of the following

\begin{prop}\label{prop5.1}
Fix some $M>1$. There exists $\La\in\ca_2\sm\ca_1$ such that $i(\La)=M$.
\end{prop}
\begin{Proof}
We construct inductively a countable family $\{\mathfrak{D}_n\}$, $n=1,2,\ld$ of sets $\mathfrak{D}_n \subset
[1,+\infty )$ according to the following rules.
\begin{enumerate}
\item[(i)] $\mathfrak{D}_1=\{ 1\}$.
\item[(ii)] $\mathfrak{D}_n=\{(a_n+\nu)^2|\,\nu=0,1,\ld,\,[a_n]+1\} \,\,\,n=1,2,\ldots$.
\item[(iii)] $\min\mathfrak{D}_{n+1}=M\cdot\max\mathfrak{D_n}$ for each $n=1,2,\ldots $,
\end{enumerate}
where $a^2_n=:\min\mathfrak{D}_n$, $n=1,2,\ld$ and $[x]$ denotes the integer part of the real number $x$ as
usual. Observe that every $n$, $m\in\N$, $n\neq m$, $\mathfrak{D}_n\cap\mathfrak{D_m}=\emptyset$. Set
$$\widetilde{\La}=\bigcup\limits^{+\infty}_{n=1}\mathfrak{D}_n.$$ We define the sequence $\La=(\la_n)$ to be the enumeration of
$\widetilde{\La}$ by the natural order.

It is obvious that $\la_n\neq0$ $\fa\;n\in\N$, $\dis\lim_{n\ra+\infty}\la_n=+\infty$, and $(\la_n)$ is a
strictly increasing sequence of positive numbers. We prove the following claim \vspace*{0.2cm} \\
\noindent
{\bf Claim 1.} For every subsequence $(\mi_n)$ of $\La$ we have
$\underset{n\ra+\infty}{\lim\sup}\Big|\dfrac{\mi_{n+1}}{\mi_n}\Big|\ge M$.
\begin{Proof}
Let some fixed subsequence $(\mi_n)$ of $\La$. Firstly we prove that for every natural number $m\in\N$, there
exists $N\in\N$ with $N\ge m$ such that $\Big|\dfrac{\mi_{N+1}}{\mi_N}\Big|\ge M$. So, fix $m\in\N$. Let $m_1$
be the unique positive integer such that $\mi_m\in\mathfrak{D}_{m_1}$. We set
$A_{m_1}:=\{n\in\N|\mi_n\in\mathfrak{D}_{m_1}\}$. It is obvious that $A_{m_1}\neq\emptyset$, since $m\in
A_{m_1}$. We set $m_2:=\max A_{m_1}$. Then $\mi_{m_2+1}\notin\mathfrak{D}_{m_1}$ and so
$\mi_{m_2+1}\ge\min\mathfrak{D}_{m_1+1}$. We have $\mi_{m_2}\le\max\mathfrak{D}_{m_1}$. Thus,
\[
\frac{\mi_{m_2+1}}{\mi_{m_2}}\ge\frac{\min\mathfrak{D}_{m_1+1}}{\max\mathfrak{D}_{m_1}}= M \ \ \text{and} \ \
m_2\ge m_1.
\]
So we proved that for every $m\in\N$, there exists some $N\ge m$ such that
\[
\frac{\mi_{N+1}}{\mi_N}\ge M.
\]
We apply the previous to an induction argument. For $m=1$ there exists $k_1\in\N$, $k_1\ge1$ such that
$\dfrac{\mi_{k_1+1}}{\mi_{k_1}}\ge M$. For $m=k_1+1$, there exists some $k_2\ge k_1+1$ (especially $k_2>k_1$)
such that $\dfrac{\mi_{k_2+1}}{\mi_{k_2}}\ge M$. Suppose that for some $\nu\in\N$ we have constructed some
$k_\nu\in\N$, such that $\dfrac{\mi_{k_\nu+1}}{\mi_{k_\nu}}\ge M$. Then for $m=k_\nu+1$ there exists some
$k_{\nu+1}\ge k_\nu+1$ (especially $k_{\nu+1}>k_\nu$) such that $\dfrac{\mi_{k_{\nu+1}+1}} {\mi_{k_{\nu+1}}}\ge
M$. Thus, we constructed a subsequence $(\mi_{k_\nu})$, $\nu=1,2,\ld$ of $(\mi_n)$ such that $k_{\nu+1}>k_\nu$
for each $\nu=1,2,\ld$ and $\dfrac{\mi_{k_\nu+1}}{\mi_{k_\nu}}\ge M$. This gives
$\dis\lim_{\nu\ra+\infty}\dfrac{\mi_{k_\nu+1}}{\mi_{k_\nu}}\ge M$, which in turn implies
$\underset{n\ra+\infty}{\dis\lim\sup}\dfrac{\mi_{n+1}}{\mi_n}\ge M$.
\end{Proof}
We show now the following \vspace*{0.2cm} \\
\noindent {\bf Claim 2.} $\underset{n\ra+\infty}{\lim\sup}\dfrac{\la_{n+1}}{\la_n}=M$.
\begin{Proof}
First of all we prove that $\dis\lim_{n\ra+\infty}a_n=+\infty$ where $a^2_n=\min\mathfrak{D}_n$, $n=1,2,\ld\;.$
Let some $n\in\N$, $n\ge2$. We have
\[
a_{n+1}=\sqrt{M}(a_n+[a_n]+1)>2a_n\sqrt{M}>2a_n\;\Rightarrow\;\frac{a_{n+1}}{a_n}>2 \ \ \text{where} \ \
a_2=\sqrt{M}>1.
\]
This gives that $\dis\lim_{n\ra+\infty}a_n=+\infty$.\medskip

Let some fixed $n\in\N$. If there exists some $m\in\N$ such that $\la_n$, $\la_{n+1}\in\mathfrak{D}_m$ then by
the construction of $\mathfrak{D}_m$ we have $\la_n=(a_m+k)^2$, $\la_{n+1}=(a_m+n+1)^2$ for some
$k\in\{0,1,\ld,[a_m]+1\}$; thus
\begin{equation} \label{sec5eq1}
\frac{\la_{n+1}}{\la_n}=\frac{(a_m+n+1)^2}{(a_m+n)^2}=\bigg(1+\frac{1}{a_m+n}\bigg)^2<
\bigg(1+\frac{1}{a_m}\bigg)^2.
\end{equation}
If there exists no $m\in\N$ such that $\la_n,\la_{n+1}\in\mathfrak{D}_m$, then this happens only if
$\la_n=\max\mathfrak{D}_m$ and $\la_{n+1}=\min\mathfrak{D}_{m+1}$ for some $m\in\N$. In this case we have
\begin{equation} \label{sec5eq2}
\dfrac{\la_{n+1}}{\la_n}=M.
\end{equation}
By (\ref{sec5eq1}), (\ref{sec5eq2}) and since $\dis\lim_{n\ra+\infty}a_n=+\infty$ the conclusion follows. This
completes the proof of Claim 2.
\end{Proof}
Claims 1 and 2 imply that $i(\La)=M$. We now show the following claim.\vspace*{0.2cm} \\
\noindent
{\bf Claim 3.} $\ssum^{+\infty}_{n=1}\dfrac{1}{\la_n}<+\infty$.
\begin{Proof}
Let some $m\in\N$, $m\ge2$. Observe that
\[
\mathfrak{D}_m=\{(a_m+\nu)^2|\,\nu=0,1,\ld,\,[a_m]+1\}.
\]
Set
\[
S_m:=\sum^{[a_m]+1}_{\nu=0}\frac{1}{(a_m+\nu)^2}.
\]
We have $a_m>1$ and for $\nu_0\in\N$
\begin{align*}
\sum^{\nu_0}_{\nu=0}\frac{1}{(a_m+\nu)^2}<\sum^{\nu_0}_{\nu=0}\frac{1}
{(a_m-1+\nu)(a_m+\nu)}&=\sum^{\nu_0}_{\nu=0}\bigg(\frac{1}{a_m-1+\nu}-
\frac{1}{a_m+\nu}\bigg) \\
&=\frac{1}{a_m-1}-\frac{1}{a_m+\nu_0}.
\end{align*}
So,
\begin{align*}
\sum^{+\infty}_{\nu=0}\frac{1}{(a_m+\nu)^2}&=\lim_{k\ra+\infty}
\bigg(\sum^k_{\nu=0}\frac{1}{(a_m+\nu)^2}\bigg)\\
&\le\lim_{k\ra+\infty}\bigg(\frac{1}{a_m-1}-\frac{1}{a_m+k}\bigg)=\frac{1}{a_m-1}.
\end{align*}
This gives
\[
S_m<\frac{1}{a_m-1}.
\]
In the proof of Claim 2 we showed that $a_{n+1}>2a_n$ for every $n\ge2$. Thus
$a_{n+1}-1>2a_n-2=2(a_n-1)\Rightarrow\dfrac{1}{a_{n+1}-1}<\dfrac{1}{2(a_n-1)}$ for $n\ge2$. The previous
inequality and an easy induction argument imply
\[
\frac{1}{a_{m+k}-1}<\frac{1}{2^k(a_m-1)},\,\,\, k=1,2,\ld
\]
Hence,
\begin{align*}
\sum^k_{\nu=0}S_{m+\nu}&=S_m+\sum^k_{\nu=1}S_{m+\nu}\\
&<S_m+\sum^k_{\nu=1}\frac{1}{a_{m+\nu}-1}\\
&<S_m+\sum^k_{\nu=1}\frac{2}{2^\nu(a_m-1)} \\
&<S_m+\frac{1}{a_m-1}\cdot\sum^{+\infty}_{\nu=1}\frac{1}{2^\nu} \\
&=S_m+\frac{1}{a_m-1}.
\end{align*}
for every $k=1,2,\ld$.  Since the sequence $\Big(\ssum^k_{\nu=0}S_{m+\nu}\Big)_k$ is strictly increasing and
bounded above by $S_m+\frac{1}{a_m-1}$ we conclude that
\[
\sum^{+\infty}_{n=m}\frac{1}{\la_n}\le S_m+\frac{1}{a_m-2}.
\]
This completes the proof of Claim 3.
\end{Proof}
The above claim shows that $\La\notin\ca_1$. We now show our last \vspace*{0.2cm} \\
\noindent
{\bf Claim 4.} $\La\in\ca_2$.
\begin{Proof}
%
%
%
%
%

For every $n\in\N$, $n\ge2$ we have
\begin{align*}
&\frac{1}{a^2_n}+\frac{1}{(a_n+1)^2}+\cdots+\frac{1}{(a_n+[a_n]+1)^2}
>\frac{1}{a_n(a_n+1)}+\frac{1}{(a_n+1)(a_n+2)}\\
&+\cdots+ \frac{1}{(a_n+[a_n]+1)(a_n+[a_n]+2)}=\sum^{[a_n]+1}_{k=0}\bigg(\frac{1}
{a_\nu+k}-\frac{1}{a_\nu+k+1}\bigg) \\
&=\frac{1}{a_n}-\frac{1}{a_n+[a_n]+2}=\frac{[a_n]+2}{a_n(a_n+[a_n]+2)}
>\frac{1}{a_n+[a_n]+2}\ge\frac{1}{2(a_n+1)}>\frac{1}{4a_n},
\end{align*}
thus,
\begin{equation} \label{sec5eq5}
a^2_n\cdot\sum^{[a_n]+1}_{k=0}\frac{1}{(a_n+k)^2}>\frac{a_n}{4}.
\end{equation}
Take now any $n\in \mathbb{N}$, $n\geq 2$. Then $\la_n\in\mathfrak{D}_m$ for some $m\in\N$, $m\ge2$. We have
\[
\la_n=(a_m+k)^2 \ \ \text{for some} \ \ k\in\{0,1,\ld[a_m]+1\}.
\]
Observe that $\la_\rho=a_{m+1}^2:=\min\mathfrak{D}_{m+1}$ for some positive integer $\rho$. Then,
\begin{align*}
\la_n\cdot\sum^{+\infty}_{k=n}\frac{1}{\la_k}&=\la_n\bigg(\sum^{[a_m]+1}_{j=k}
\frac{1}{(a_m+j)^2}+\sum^{+\infty}_{\nu=\rho}\frac{1}{\la_\nu}\bigg) \\
&=\la_n\cdot\sum^{[a_m]+1}_{j=k}\frac{1}{(a_m+j)^2}+\la_n\bigg(
\sum^{+\infty}_{\nu=\rho}\frac{1}{\la_\nu}\bigg) \\
&=\la_n\sum^{[a_m]+1}_{j=k}\frac{1}{(a_m+j)^2}+\frac{\la_n}{a_{m+1}^2}
\bigg(\la_\rho\cdot\sum^{+\infty}_{\nu=\rho}\frac{1}{\la_\nu}\bigg)\\
&>\frac{\la_n}{a_{m+1}^2}\cdot\bigg(\la_\rho\cdot\sum^{+\infty}_{\nu=\rho}\frac{1}{\la_\nu}
\bigg)>\frac{1}{9M}\cdot\frac{a_{m+1}}{4},
\end{align*}
where the last inequality above follows by (\ref{sec5eq5}) and the inequality $9Ma_m^2>a_{m+1}^2$ which is easy
to show and it is left to be checked by the reader. The previous estimate and the fact that
$\dis\lim_{m\ra+\infty}a_m=+\infty$ implies the desired result. This completes the proof of Claim 4. Hence
$\La\in\ca_2$ and the proof of this proposition is complete. \qb
\end{Proof}
\end{Proof}

{\bf Acknowledgements}. I am grateful to George Costakis for his helpful comments and remarks and for all the
help he offered me concerning the presentation of this work.

\end{document}